[1994/12/01]
\documentclass[twoside,reqno,final,11pt]{amsart}
\usepackage{amsmath,amssymb,amsthm,amscd}
\usepackage{graphicx}

\def\arxivno#1{\href{http://arxiv.org/#1}{\texttt{#1}}}

\usepackage[
    pdfauthor={Sergei Yakovenko},
    bookmarksopen=false,
    pdfstartview=FitH,
    linkcolor=blue,
    citecolor=blue,
    urlcolor=blue,
     plainpages,
    ]{hyperref}

\numberwithin{equation}{section}

\newtheorem{Thm}{Theorem}
\newtheorem{Lem}{Lemma}
\newtheorem{Cor}[Lem]{Corollary}
\newtheorem{Prop}[Lem]{Proposition}

\theoremstyle{definition}

\newtheorem{Ex}{Example}

\theoremstyle{remark}
\newtheorem{Rem}{Remark}

\let\parasymbol=\S
\def\secref#1{\parasymbol\ref{#1}}
\def\S{\varSigma}
\def\e{\varepsilon}
\def\ssm{\smallsetminus}

\renewcommand\ge{\geqslant}
\renewcommand\le{\leqslant}
\let\tildeaccent=\~
\let\hataccent=\^
\renewcommand\~[1]{\widetilde{#1}}
\renewcommand\^[1]{\widehat{#1}}

\def\<{\left<}
\def\>{\right>}
\def\({\ifmmode\left(\else\textup{(}\fi}
\def\){\ifmmode\right)\else\textup{)}\fi}

\def\C{\mathbb C}
\def\R{\mathbb R}
\def\Mat{\operatorname{Mat}}
\def\:{\colon}
\def\l{\lambda}
\def\pd#1#2{\frac{\partial #1}{\partial #2}}
\def\Re{\operatorname{Re}}
\def\Im{\operatorname{Im}}

\begin{document}
\DeclareGraphicsExtensions{.jpg,.pdf,.mps,.png,.gif,.tif}

\title [Grigoriev Theorem]
 { Oscillation of linear
 ordinary differential equations: \\ on a theorem
 by A.~Grigoriev }

\author[S.~Yakovenko]
 { Sergei Yakovenko$^*$ }

\thanks{$^*$Gershon Kekst professorial chair in Mathematics.}


\address{\href{http://www.wisdom.weizmann.ac.il}
 {Weizmann Institute of Science}\\
 Rehovot\\
 Israel}
\email{\href{mailto:yakov@wisdom.weizmann.ac.il}
 {\tt yakov@wisdom.weizmann.ac.il}
 \endgraf{\leavevmode\kern-2pt\it WWW page\/}:
 \href{http://www.wisdom.weizmann.ac.il/~yakov}
 {\tt http://www.wisdom.weizmann.ac.il/\char'176 yakov}}

 \date{September 2004}


\begin{abstract}
 We give a simplified proof and an improvement of a recent
 theorem by A.~Grigoriev, placing an upper bound for the number of
 roots of linear combinations of solutions to systems of linear
 equations with polynomial or rational coefficients.
\end{abstract}

 \maketitle

\section{Background on counting zeros of solutions
of linear ordinary differential equations}

\subsection{De la Val\'ee Poussin theorem and Novikov's
counterexample} A linear $n$th order homogeneous differential
equation
\begin{equation}\label{lode}
    y^{(n)}+a_1(t)\,y^{(n-1)}+\cdots+a_{n-1}(t)\,y'+a_n(t)\,y=0,
    \qquad y^{(k)}=\tfrac{d^ky}{dt^k},
\end{equation}
with real continuous coefficients $a_j(t)$ is called
\emph{disconjugate} (\emph{Chebyshev, non-oscillating}) on a real
interval $[\alpha,\beta]\subseteq\R$, if any nontrivial solution
$y(t)$ of this equation has at most $n-1$ isolated root on this
interval.

A theorem by C.~de la Vall\'ee Poussin \cite{poussin} asserts that
any equation \eqref{lode} is disconjugate on any interval
sufficiently short relative to the \emph{magnitude} of the
coefficients of the equation. More precisely, if $|a_j(t)|\le b_j$
on $[\alpha,\beta]$ and $\sum_{j=1}^n b_j\,|\beta-\alpha|^j/j!<1$,
then \eqref{lode} is disconjugate. This allows to place a rather
accurate upper bound on the number of isolated roots of any
solution of a known differential equation in terms of the length
of the interval on which the solution is considered and the
magnitude of the coefficients of the equation. A complex analog of
this theorem was obtained in \cite{fields} for linear homogeneous
$n$th order equations with holomorphic coefficients in a polygonal
complex domain $t\in U\Subset\C$.

A natural generalization of these results for \emph{systems} of
first order linear homogeneous differential equations of the form
\begin{equation}\label{ls-1}
    \dot x(t)=A(t)\,x(t),\qquad
    x=\begin{pmatrix}x_1\\ \vdots \\ x_n\end{pmatrix},
    \quad
    A(t)=\begin{pmatrix}a_{11}(t)&\cdots&a_{1n}(t)
    \\
    \vdots&\ddots&\vdots
    \\
    a_{n1}(t)&\cdots&a_{nn}(t)
    \end{pmatrix},
\end{equation}
with variable coefficients matrix $A(t)$, might deal with a
uniform upper bound for the number of isolated zeros of an
arbitrary linear combination of components
$y(t)=c_1\,x_1(t)+\cdots+c_n\,x_n(t)$, in terms of the length of
the interval $t\in[\alpha,\beta]$ where the solution is
considered, and the magnitude $\max_{t\in[\alpha,\beta]}|A(t)|$ of
the matrix norm on this interval. Among other reasons, such result
could be expected since for any such linear combination one can
explicitly derive a linear $n$th order equation of the form
\eqref{lode}, as explained in \secref{sec:sys2eq}.

Unfortunately, nothing like this is possible. In
\cite{mit:counterexample} D.~Novikov constructed a very simple
system with the following properties:

\begin{enumerate}
\item the coefficient matrix $A(t)$ is a real matrix polynomial
of degree $m$,
\begin{equation}\label{ls-p}
    A(t)=\sum_{j=0}^{m-1}A_k\,t^k,\qquad A_k\in\Mat_n(\C),
    \quad k=0,\dots,m-1,
\end{equation}
in the minimal dimension $n=2$,

\item the matrix norm $|A(t)|\le 1$ everywhere in the disk $\{|t|<2\}\subset\C$;

\item the first component $x_1(t)$ of
some solution $x(t)$ of the respective system \eqref{ls-1} has as
many isolated roots in the unit disk $\{|t|<1\}$ as necessary,
provided that the degree $m$ of the matrix polynomial can be
arbitrarily large.
\end{enumerate}

When reducing the Novikov's system to a second order equation
\eqref{lode}, one encounters uncontrollably large rational
coefficients $a_j(t)=p_j(t)/p_0(t)$: while the numerators $p_j(t)$
can be explicitly majorized, the denominator $p_0(t)$ can be
arbitrarily small (either coefficients-wise, or uniformly on large
intervals).

\subsection{Bounded meandering principle}
In an attempt to circumvent the appearance of the small
denominators, one can declare all entries of the matrix
coefficients $A_0,A_1,\dots,A_{m-1}$ of the  polynomial
\eqref{ls-p} to be new \emph{phase variables} $\l\in\R^{n^2m}$
playing the same role as the initial variables $x\in\R^n$. Then
the system \eqref{ls-p} coupled with the trivial equations
$\dot\l=0$, $\dot t=1$, becomes \emph{nonlinear}, but all numeric
coefficients of this nonlinear system become \emph{integer} (in
fact, only zeros or ones). This \emph{lattice structure} allows to
avoid appearance of the small denominators (any integer number is
either zero or at least one in the absolute value), but the
``reduction'' to a scalar \emph{linear} equation becomes
considerably more complicated.

In \cite{annalif-99} the problem of counting roots of solutions of
systems of polynomial differential equations was reduced to
estimating the length of an ascending chain of polynomial ideals.
The ultimate result established in \cite{annalif-99} is a general
``bounded meandering principle'': \emph{for any system of
polynomial ordinary differential equations and any polynomial
combination of components of its solutions, the number of isolated
roots of this combination on any finite interval can be explicitly
bounded in terms the degrees of all polynomials, dimension of the
system, the length of the interval and the magnitude of the
coefficients of the equations}. More specifically, if the
coefficients of the system of $n$ differential equations and the
length of the interval on which the combination is considered, do
not exceed some upper bound $M>2$, and the degrees of all
polynomials do not exceed $m$, then the number of isolated roots
of any polynomial combination of solutions does not exceed $M^q$,
where $q=q(n,m)$ is an explicit computable function of two natural
arguments, growing no faster than the iterated exponent of order
$4$:
\begin{equation}\label{tower4}
    q(n,m)\le 2^{2^{2^{2^{4n\ln m+O(1)}}}}.
\end{equation}
where $O(1)$ is an absolute constant \cite{annalif-99}.

\subsection{Grigoriev theorem}
The approach based on the bounded meandering principle, completely
ignores the linearity of the initial system \eqref{ls-p}: the
bound for the number of isolated zeros is given by the same
extremely excessive expression \eqref{tower4}.

In 2001 in his Ph.~D.~thesis \cite{alexg:thesis} recently
published in the revised form \cite{alexg:arxiv}, Alexei Grigoriev
returned to the study of the system \eqref{ls-p} by ``purely
linear'' methods. The principal obstruction on this way, already
mentioned above, is occurrence of the ``singular degeneracy'': for
certain values of the parameters $\l$, the resulting equation has
the leading coefficient (before the highest derivative) vanishing
\emph{identically}. In the classical language, the parametric
family of differential equations near such values of the
parameters would exhibit a \emph{singular perturbation}.

Appearance of singular perturbations can indeed be a source of
highly oscillating solutions. The simplest examples can be easily
constructed already in the class of linear equations with constant
coefficients. Yet in the situation when the equation is obtained
from a system \eqref{ls-p}, we have an additional bit of
information: \emph{all solutions of the derived equation depend
analytically on the parameters}.

In \cite[p.~520]{fields} it was conjectured that even in the
family exhibiting a singular perturbation, all solutions which
depend analytically on the parameters admit an effective
\emph{uniform} upper bound for the number of isolated zeros.
A.~Grigoriev proved a similar conjecture for the equations with
rational coefficients whose all solutions remain analytic, and
achieved a much better upper bound than the tower \eqref{tower4}.
In the simplest settings his result can be formulated for a system
of the form \eqref{ls-p} with the \emph{polynomial} coefficients.

Being polynomial, such system can be considered also for complex
values of time $t$. We are interested in an upper bound for the
number of isolated roots of any linear combination in a bounded
domain. Without loss of generality we may assume (rescaling the
variable $t$, if necessary) that this domain is the unit disk
$\mathbb D=\{|t|<1\}$, while all matrix coefficients $A_k$ are
bounded by the constant $M$:
\begin{equation}\label{ls-r}
    |\,t\,|< 1,\qquad |A_k|\le M,\quad k=0,\dots,m-1.
\end{equation}
The following result can be obtained by slightly improving the
arguments by A.~Grigoriev.

\begin{Thm}\label{main}
The number of isolated zeros of any linear combination of
components of solutions of the system \eqref{ls-p} constrained by
the restrictions \eqref{ls-r} in the unit disk does not exceed
\begin{equation}\label{main-bd}
        [\max(M,2)]^{2^{O(n^2m)}},
\end{equation}
where the constant in $O$ is absolute.
\end{Thm}

\begin{Rem}
The bound \eqref{main-bd} seems to be still rather excessive. For
instance, using the standard arguments relating the growth of
analytic functions with the number of their zeros via the Jensen
inequality \cite{levin,jde-96}, we can prove that \emph{for most}
linear combinations, the number of zeros does not exceed $M\cdot
2^{O(m)}$, if we consider $n$ as a fixed parameter. This is a far
shot from the double exponential bound \eqref{main-bd}.
\end{Rem}

\begin{Rem}
Grigoriev in \cite{alexg:arxiv} formulates this theorem only for
the case $M=1$, and only as a remark; besides, the bound that he
claims in this case, is somewhat worse: $2^{2^{(nm)^{O(1)}}}$
instead of $2^{2^{O(n^2m)}}$ in \eqref{main-bd}.
\end{Rem}

The goal of this work is to reveal some geometric ideas lying
behind the Grigoriev's proof which (despite the recent revision)
remains still highly technical and sometimes rather obscure.
Besides, we discuss how some other classes of differential
equations can be treated in the similar way, in particular, why
the hypergeometric systems are better than a general Fuchsian
system from the point of view of counting the roots. This is
especially important in connection with the tangential Hilbert
problem on the number of zeros of Abelian integrals
\cite{redundant,montreal}.

\section{From system to a high order equation}\label{sec:sys2eq}

\subsection{Variations on the theme of de la Vall\'ee Poussin}
In what follows we will use the $\ell^1$-norm on the space of
polynomials in one or several variables: by definition, the norm
of a polynomial is the sum of absolute values of all (real or
complex) coefficients,
\begin{equation}\label{el1-def}
    f\in\C[x]=\C[x_1,\dots,x_n],\quad f=\sum_{|\alpha|\le d}c_\alpha
    x^\alpha,\qquad \|f\|_{\C[x]}=\sum_{\alpha}|c_\alpha|.
\end{equation}
The main advantage of this norm is its multiplicativity:
$\|fg\|\le\|f\|\,\|g\|$.

Consider the equation \eqref{lode} with \emph{polynomial}
coefficients $a_0,a_1,\dots,a_n\in\C[t]$ of known degrees. Note
that, since the equation is homogeneous, the coefficients are
defined modulo a common factor. The following result is a
variation on the theme of Theorem 2$'$ from \cite{jde-96}: if the
principal coefficient before the highest derivative is not too
\emph{uniformly} small while all other coefficients explicitly
bounded from above, then the number of zeros of any
\emph{holomorphic} solution is explicitly bounded (non-holomorphic
solutions with singularities may disobey this principle).

\begin{Lem}\label{lem:ode}
Assume that the coefficients of the equation \eqref{lode} of order
$n$ are polynomials of degree $\le d$, and their norms
$\|a_j\|=\|a_j\|_{\C[t]}$ satisfy the restrictions
\begin{equation}\label{norms}
    \|a_0\|=1,\qquad \|a_k\|\le K,\quad k=1,\dots,n.
\end{equation}

Then any \emph{entire} solution to \eqref{lode} may have at most
$Kn\cdot 2^{O(d)}$ isolated roots in the unit disk $\mathbb
D=\{|t|<1\}\subset\C$, where the constant in $O(d)$ is absolute.
\end{Lem}

\begin{proof}
We start with observing that a polynomial $a_0\in\C[t]$ of degree
$\le d$ normalized by the condition $\|a_0\|=1$, cannot be
uniformly small on, say, the annulus $\{1<|t|<2\}$.  More
precisely, there exists a circle $\{|t|=r\}$, $1<r<2$, on which
$|a_0(t)|$ is bounded from \emph{below}, $|a_0(t)|\ge 2^{-O(d)}$.

Indeed, the condition $\|a_0\|=1$ means that one of the
coefficients of the polynomial $a_0$ is at least $1/d$ in the
absolute value. Applying the Cauchy formula for the respective
derivative of $a_0$, we conclude that
\begin{equation*}
    1/d\le\left|\frac1{2\pi
    i}\oint_{|s|=1}\frac{a_0(s)}{s^{k+1}}\,ds\right|
    \le \max_{|s|=1}|a_0(s)|.
\end{equation*}
Therefore that there exists a point $s_0$ on the unit circle
$|t|=1$, on which $|a_0(s_0)|$ is at least $1/d$. Applying
\cite[Theorem 4, p.~79]{levin} to the polynomial
$\~a_0(s)=a_0(s-s_0)/d$, we conclude that the lower bound
$|a_0(s)|\ge 2^{-O(d)}$ with some absolute constant in $O(d)$
holds on the complement to a finite union of disks of the total
diameter less than $1$. Such disks cannot intersect all concentric
circles $\{|t|=r\}$, $1<r<2$, which proves the required assertion.

On this circle $\{|t|=r\}$ all other coefficients are bounded from
above, $|a_j(t)|\le K\cdot 2^{O(d)}$ for all $j=1,\dots,n$.
Restricting the equation \eqref{lode} on the circle parameterized
by $t=r\exp is$, $s\in[0,2\pi]$, and applying Corollary~2.7 from
\cite[p.~507]{fields}, we conclude that the variation of argument
of \emph{any} solution $y(t)$ of the equation \eqref{lode} on this
circular segment, can be at most $Kn\cdot 2^{O(d)}$.

If $y(t)$ is holomorphic in the disk $\{|t|<r\}$, in particular,
if this solution is an entire function, then by the argument
principle the number of its zeros does in the smaller disk
$\mathbb D$ does not exceed $Kn\cdot2^{O(d)}$, as asserted.
\end{proof}

\begin{Rem}\label{rem:larger-disk}
The requirement that the solution is entire, is excessive. In
fact, it is sufficient that the solution of the equation is
analytic in a disk of any finite radius larger than $1$. For
instance, the proof below works literally if $y(\cdot)$ is
analytic in the disk of radius $6$. If necessary, the number $6$
can be replaced with only minor modifications by any number bigger
than $1$, but this choice will affect the constant in $O(d)$.
\end{Rem}

\subsection{Derivation of the principal equation}\label{sec:deriv}
Consider now the linear system \eqref{ls-p}. We will denote by
$\l=\{\l_{ijk}\}\in\C^{n^2m}$ the collection of all entries of all
matrices $A_k$, so that the coefficients matrix of the system
\eqref{ls-p} would be a ``universal'' matrix polynomial $A(t,\l)$
from the ring $\mathbb Z[t,\l]$, whose entries have norm $m$ in
this ring.

Without loss of generality, when proving Theorem~\ref{main}, it is
sufficient to consider only the first coordinate $y(t)=x_1(t)$ of
an arbitrary solution of \eqref{ls-p}.

The formulas for higher derivatives of the function $y(\cdot)$ can
also be expressed by linear combination of the components
$x_1(t,\l),\dots,x_n(t,\l)$ with polynomial coefficients. Indeed,
assume inductively that
\begin{equation*}
    y^{(k)}(t)=\xi_k(t,\l)\cdot x(t,\l),
\end{equation*}
where $\xi_k(t,\l)$, $k=0,1,\dots,n$ is a polynomial (row) vector
function. Then after derivation by virtue of the system
\eqref{ls-1} we obtain
\begin{equation*}
    y^{(k+1)}(t)=\pd{\xi_k}t\cdot x(t,\l)+\xi_k\cdot A(t,\l)\cdot
    x(t,\l)=\(\pd{\xi_k}t+\xi_k A\)\cdot x(t,\l).
\end{equation*}
This allows to define the vector functions $\xi_k$ starting from
from the initial vector-function
$\xi_0(t,\l)\equiv(1,0,\dots,0)\in\C^n$.

Consider the collection $\xi_k\in\C[t,\l]\otimes\C^n$ of
polynomial vector-functions defined by the recursive formula
\begin{equation}\label{recurs}
    \xi_{k+1}=\pd{}t \xi_k+\xi_k\cdot A(t,\l),\qquad
    k=0,1,\dots,n-1.
\end{equation}

The process \eqref{recurs} is so explicit that the following
estimates are almost obvious.

\begin{Lem}\label{iter}
The collection of polynomial vector functions
$\xi_0,\dots,\xi_n\in\mathbb Z[t,\l]\otimes\C^n$ satisfies the
following inequalities\textup:
\begin{enumerate}
 \item the vector polynomials $\xi_0,\dots,\xi_n$ are integral,
 i.e., they belong to the lattice $\mathbb Z[t,\l]\otimes_{\mathbb Z}
 \mathbb Z^n\subset\C[t,\l]\otimes_\C\C^n$,
 \item $\deg\xi_k\le km$ for all $k=0,\dots,n$,
 \item the wedge product of any $n$ vectors out of this collection
is a scalar integral polynomial from $\mathbb Z[t,\l]$ of the
degree at most $\tfrac12n(n-1)m$ and the norm not exceeding
$n!\,(2nm)^{n^2}$,
 \item the wedge product $\xi_0\land\cdots\land\xi_{n-1}\in\mathbb
Z[t,\l]$ is not identically zero in $\mathbb Z[t,\l]$.
\end{enumerate}
\end{Lem}

\begin{proof}
Everywhere until the end of the proof we abbreviate
$\deg=\deg_{t,\l}$ and $\|\cdot\|=\|\cdot\|_{\C[t,\l]}$.

The fact that all coefficients of $\xi_k$ are \emph{integral} over
$t,\l$, follows immediately from the fact that the ``universal''
matrix polynomial $A(t,\l)$ is integral, $A\in\mathbb
Z[t,\l]\otimes_{\mathbb Z}\Mat(n,\mathbb Z)$.

The inequality
\begin{equation*}
    \deg\xi_{k+1}\le\deg\xi_k+\deg A(t,\l)\le
    \deg\xi_k+m
\end{equation*}
proves the assertion about the degrees of the functions $\xi_k$
themselves. The Laplace formula for the complete expansion of the
$(n\times n)$-determinant proves the assertion for the degrees of
the wedge products.

The assertion about the norm follows from the recursive inequality
\begin{equation*}
    \|\xi_{k+1}\|\le(\deg\xi_k+nm)\|\xi_k\|\le 2nm\|\xi_k\|,
\end{equation*}
where by the norm of the polynomial \emph{vector}-function we
understand the maximum of the norms of its scalar polynomial
components. Therefore each entry of the $(n\times n)$-determinant
has the norm not exceeding $(2nm)^n$, and by the Laplace expansion
we have the sum of $n!$ terms, each of the norm at most
$[(2nm)^n]^n$ by the multiplicativity of the norm.

To prove the last assertion on the nondegeneracy, it is sufficient
to note that for the specific value of the variables $\l$
corresponding to the ``\emph{constant}'' (i.e., independent of
$t$) polynomial $A=A_0$ without higher terms, the vectors
$\xi_k=A_0^k\xi_0$ are also ``constant'' vectors. For the generic
choice of $A_0\in\Mat(n,\C)$ they are obviously linear independent
over $\C$ for all $k=0,1,\dots,n-1$, and hence the wedge product
is not identically zero.
\end{proof}

\begin{Cor}\label{eq-der}
The first coordinate of any solution of the system \eqref{ls-p} is
an entire function satisfying a linear $n$th order differential
equation of the form
\begin{equation}\label{meq}
    a_0(t,\l)\,y^{(n)}+a_{1}(t,\l)\,y^{(n-1)}+\cdots+a_n(t,\l)\,y=0,
\end{equation}
for the corresponding value of the parameters $\l$.

The coefficients $a_j(t,\l)\in\mathbb Z[t,\l]$ of this equation
satisfy the following conditions\textup:
\begin{enumerate}
 \item all $a_j\in\mathbb Z[t,\l]$ are \emph{integral} \textup(lattice\textup)
 polynomials,
 \item $\deg_{t,\l}a_j\le\tfrac12n(n-1)m$,
 \item $\|a_j\|_{\C[t,\l]}\le n!\,(2nm)^{n^2}$,
 \item the leading coefficient $a_0(t,\l)$ does not
vanish identically.
\end{enumerate}
\end{Cor}

\begin{proof}
For arguments of pure dimensionality, the vectors
$\xi_0,\dots,\xi_n\in\mathbb Z[t,\l]\otimes\mathbb
Z^n\subset\C(t,\l)\otimes\C$ must be linear dependent over the
field of rational functions $\C(t,\l)$,
\begin{equation*}
    a_0\xi_n+\cdots+a_n\xi_0=0.
\end{equation*}
This identity implies a linear equation of order $n$ for the
function $y(t)=\xi_0\cdot x$. The coefficients of this dependence
can be found (by the Cramer rule) as the ratios of the wedge
products of the vectors $\xi_j$ of degree $n$. These were already
described in Lemma~\ref{iter}.
\end{proof}

\subsection{Investigation of the principal equation}
The fact that the leading coefficient $a_0$ of the equation
\eqref{meq} does not vanish \emph{identically} in $\mathbb
Z[t,\l]$, \emph{does not imply} that for some values of the
\emph{parameters} $\l$ it becomes a zero polynomial of $t$. Such
values of the parameters correspond to degenerated systems for
which the linear dependence over the field $\C(t)$ between
derivatives of the first coordinate occurs earlier than on the
$n$th step. In classical terminology, this phenomenon should be
described as a \emph{singular perturbation}, when the leading
coefficient of a differential equation vanishes identically.

The main result by Grigoriev, which was conjectured in
\cite[\parasymbol 4.5, p.~520]{fields}, is the assertion that in
fact this singular perturbation is only \emph{apparent}, namely,
that for these exceptional values of the parameter, all other
coefficients of the principal equation \eqref{meq} are also
vanishing identically, and their ratios are explicitly bounded for
all nearby values of the parameters $\l$.

More precisely, consider the functions
\begin{equation}\label{nm}
    b_j\:\C^{n^2m}\to\R_+,
    \qquad
    b_j(\l)=\|a_j(\cdot,\l)\|_{\C[t]},\qquad j=0,1,\dots,n,
\end{equation}
the norms of the coefficients $a_j(t,\l)$ of the equation
\eqref{meq} considered as \emph{univariate} polynomials. The zero
locus $\S$ of the function $b_0$ corresponds to the ``singular
degeneracy'' of this equation.

The cornerstone of the proof of Theorem~\ref{main} is the
following analytic result discussed in details and proved in the
next section \secref{sec:1-param}.

\begin{Lem}[Principal Lemma]\label{grig}
The ratios $b_j/b_0$ are locally bounded on $\C^{n^2m}\ssm\S$ for
all $j=1,\dots,n$.
\end{Lem}

In the remaining part of this section we derive Theorem~\ref{main}
from Lemma~\ref{grig}. The first step is the following Corollary,
transforming the quantitative statement of Lemma~\ref{grig} into a
qualitative bound.

\begin{Lem}\label{cor:bound}
In the notations \eqref{meq}--\eqref{nm},
\begin{equation}\label{bds}
    \sup_{|\l|\le M}\frac{b_j(\l)}{b_0(\l)}\le
    \max(M,2)^{2^{O(n^2m)}},\qquad j=1,\dots,n,
    \end{equation}
where the supremum is taken over all values
$\l\notin\S=\{b_0(\l)=0\}$.
\end{Lem}

\begin{proof}[Proof of the Lemma]
Consider instead of the norms $b_j(\l)=\|a_j(\cdot,\l)\|_{\C[t]}$
equal to the sum of modules coefficients of the univariate
polynomials $a_j(\cdot,\l)$, the functions $c_j(\l)$ equal to the
sum of the squares of absolute of those coefficients. The
advantage is that the functions $c_j(\l)$ are \emph{polynomials}
in $\Re\l$ and $\Im\l$. On the other hand, by the Cauchy
inequality we have two-sided equivalence,
\begin{equation*}
   {c_j(\l)}\le b_j^2(\l)\le n\,{c_j(\l)},
\end{equation*}
therefore it will be sufficient to prove \eqref{bds} for the
ratios $c_j(\l)/c_0(\l)$: the additional power $1/2$ and the
factor $\sqrt n$ will be obviously absorbed in the $O$-symbol.

Without loss of generality we may assume that the number $M$ is
natural. In the real space $\R^{2n^2m}\times\R_+$ (the first
factor is the realification of $\C^{n^2m}$) consider the subset
\begin{equation*}
    S_{j,M}=\{(\l,u)\:c_0(\l)\ne0,\ u\ge 0,\ |\l|\le M,\ c_0(\l)\,u\le
    c_j(\l)\}.
\end{equation*}
This subset is semialgebraic (defined by polynomial equalities and
inequalities). The coefficients of these equalities and
inequalities are \emph{integer} numbers which do not exceed
$\max(M,N)$, where $N=n!(2nm)^{n^2}$, and the degrees do not
exceed $n^2m$.

In addition, this set is \emph{bounded} by Lemma~\ref{grig}:
$S_{j,M}$ is a subgraph of a locally bounded function defined on a
bounded set $\{|\l|\le M\}$.

A {bounded} semialgebraic set defined by \emph{integer} (lattice)
polynomial (in)equalities of degree $\delta$ and the norms at most
$\mu$ in $\R^{\nu}$, admits an explicit, effective and rather
accurate upper estimate of its size: by
\cite{grigoriev-vorobjov,heintz-roy-solerno,roy.e.a.}, such set
always belongs to the cube of the size $\mu^{\delta^{O(\nu)}}$
centered at the origin in $\R^\nu$.

Substituting the known values of the norms and the degrees, we
obtain the required upper bound for the $u$-coordinate, as
required. Note the most crucial parameter is the dimension $\nu$
equal in our case to $n^2m$ (the total number of different
parameters describing the ``universal'' system \eqref{ls-p}). All
other dependencies on $n,m$ are absorbed into the term $O(n^2m)$
occurring at the second floor of the tower function.
\end{proof}

\begin{Ex}
Consider the real rational function of two variables
$r(x,y)=(ax^2+by^2)/(cx^2+by^2)$. If the coefficients
$a,b,c,d\in\R$ are all positive and $ad-bc\ne0$, then this
function is locally bounded near $(x,y)=(0,0)$. Yet the limits
along different rays $y=ux$, $u\in\R$, are different and vary
between $a/c$ and $b/d$. Knowing only the upper bound $\mu$ for
the coefficients $a,b,c,d$ yields no local upper bound for
$r(x,y)$, but if these coefficients are integral, then
$\mu^{-1}\le r(x,y)\le \mu$.
\end{Ex}

\subsection{Proof of Theorem~\ref{main}}
For all values of the parameters $\l$ outside the degeneracy locus
$\S$, the assertion of Theorem~\ref{main} follows from
Lemma~\ref{lem:ode} into which one should substitute the value of
$K=\max(M,2)^{2^{O(n^2m)}}$ given by the Lemma~\ref{cor:bound}.
For $\l\in\S$ the principal equation \eqref{meq} degenerates into
the trivial identity $0=0$, but \emph{solutions} (the first
coordinates of solutions of the system \eqref{ls-p}) remain entire
functions of $t$ and $\l$. The number of isolated zeros of an
analytic function in the \emph{open} disk $\{|t|<1\}$ depends
semi-continuously on the parameters and hence the expression given
in the Theorem, serves as an upper bound for the number of
isolated zeros also for the exceptional values of the parameters
$\l\in\S$.\qed

\section{Analytic one-parameter families of linear
 subspaces}\label{sec:1-param}

The proof of the Principal Lemma~\ref{grig} is based on the
following observation: a linear space spanned by any number of
vectors analytically depending on \emph{one} parameter, depends
very regularly on this parameter even if the vectors exhibit extra
linear dependence for isolated values of the parameter.

\subsection{Finite-dimensional case}
We first illustrate the above observation for subspaces of a
finite-dimensional space. Everywhere below the parameter $\e$ is
one-dimensional (real or complex) and local, varying in a
neighborhood of zero.

\begin{Prop}\label{findim}
Let $v_1,\dots,v_m\:(\R^1,0)\to\R^n$ be any real analytic germs of
vector-functions, $m\le n$. Then there exist $k$ real analytic
germs of vector functions $w_1,\dots,w_k\:(\R^1,0)\to\R^n$, $k\le
n$, which have the constant rank $k$ and span the same linear
subspace in $\R^n$ as the vectors $v_1(\e),\dots,v_m(\e)$ for all
$\e\ne0$.
\end{Prop}

\begin{proof}
The proof goes by induction in $m$. If $m=1$ and $v_1(\e)\equiv0$,
then there is nothing to prove, otherwise $v_1(\e)=\e^\nu w_1(\e)$
with $w_1$ an analytic vector-function such that $w_1(0)\ne0$.
Clearly, the linear spaces spanned by $v_1(\e)$ and $w_1(\e)$
coincide for $\e\ne0$.

To prove the induction step, assume that the vectors
$w_1(\e),\dots,w_{k-1}(\e)$ are linear independent and span the
same linear space $L(\e)$ as the vectors
$v_1(\e),\dots,v_{m-1}(\e)$ for all $\e\ne0$. If $v_m(\e)\in
L(\e)$ for all $\e$, then $w_1,\dots,w_{k-1}$ span the same space
as $v_1,\dots,v_m$. Otherwise choose any linear subspace $K$ in
$\R^n$ complementary to $L(0)$  (hence to all $L(\e)$ for small
$\e$) and the projection $\pi_\e\:\R^n\to L(\e)$ parallel to $K$.
The vector function $v_m(\e)-\pi_\e v_m(\e)$. By assumption, it is
not identically zero and is real analytic, since $\pi_\e$ depends
analytically on $\e$ by the Implicit function theorem. Therefore,
\begin{equation*}
    v_m(\e)-\pi_\e v_m(\e)=\e^\nu w_k(\e),\qquad w_k(0)\ne 0,
\end{equation*}
for some finite power $\nu$. The vector $w_k(0)$ is transversal to
$L(0)$, hence the system $w_1(\e),\dots,w_k(\e)$ is linear
independent for all small $\e$ and generate the same space as
$v_1(\e),\dots,v_m(\e)$.
\end{proof}

Clearly, literally the same arguments prove the complex analytic
counterpart of this theorem. Besides, one may replace holomorphic
functions by meromorphic (in which case the exponents $\nu$ in the
above proof may well be negative).

\begin{Rem}
The coordinate form of Proposition~\ref{findim} may be considered
as a theorem on rank for meromorphic matrix functions (modulo
renaming the parameter to a more traditional complex variable
$z$). Any meromorphic germ of a $n\times m$-matrix function
$X(z)$, $z\in(\C^1,0)$, can be represented as
\begin{equation}\label{matr-rep}
    X(z)=U(z)D(z)V(z),\qquad
    D(z)=\begin{pmatrix}
    z^{\nu_1}&&&&\\
    &\ddots&&&\\
    &&z^{\nu_k}&&\\
    &&&0&\\
    &&&&\ddots\\
    \end{pmatrix}
\end{equation}
where the square holomorphic matrix functions $U(z)$ and $V(z)$ of
appropriate sizes are holomorphically invertible, and $D(z)$ is an
$n\times m$-matrix with only powers or zeros on the diagonal.

Another interpretation of this result can be described as an
extendability of holomorphic subbundles. Consider a holomorphic
vector bundle over a one-dimension base (curve) $Z$. Assume that a
holomorphic subbundle is defined over a punctured neighborhood of
some point $z_0\in Z$ on the base. If the subbundle has a
``moderate'' singularity, i.e., it can be spanned by (eventually,
multivalued) sections growing no faster than polynomially at
$z_0$, then in fact this subbundle can be holomorphically extended
at $z_0$.
\end{Rem}

\subsection{Infinite-dimensional version}
The above arguments need some modification to be applicable to
finite-dimensional subspaces of the infinite-dimensional space of
analytic functions of one variable. Each such subspace can be
naturally identified with a homogeneous linear ordinary
differential equation.

Everywhere below $\partial$ stands for the (partial) derivative
$\partial/\partial t$; as before, $\e\in(\C^1,0)$ is a
one-dimensional complex parameter.

\begin{Lem}\label{lode-e}
Let $f_1(t,\e),\dots,f_n(t,\e)$ be $n$ functions analytic in
$\mathbb D\times(\C^1,0)$, which are linear independent as
functions of $t\in\mathbb D$ for all $\e\ne0$.

Then the functions $f_1,\dots,f_n$ satisfy a linear $n$th order
homogeneous differential equation
\begin{equation*}
    L_\e y=0,\qquad\text{where}\quad
    L_\e=a_0(t,\e)\partial^n+\cdots+a_{n-1}(t,\e)\partial+a_n(t,\e),
\end{equation*}
with the coefficients $a_j$ analytic in $\mathbb D\times(\C^1,0)$
and \emph{non-singular} at $\e=0$, so that
$a_0(\cdot,0)\not\equiv0$.
\end{Lem}

\begin{proof}
The proof imitates that of Proposition~\ref{findim} and goes by
induction in $n$. For $n=1$ the claim is obvious: the function
$f_1(t,\e)=\e^\nu g(t,\e)$, where $g(\cdot,0)\not\equiv0$ is
analytic in $\mathbb D\times(\C^1,0)$, satisfies the linear
equation $L_1f_1=0$, where $L_1$ is the linear first order
differential operator,
\begin{equation}\label{L1}
    L_1=g\partial-(\partial g),\qquad
    y\overset{L_1}{\longmapsto}gy'-g'y.
\end{equation}
By the choice of $g$, the coefficients of this operator depend
holomorphically on $\e$ and the leading coefficient $g(t,\e)$ does
not vanish identically at $\e=0$. The base of induction is thus
established.

Assume now that the differential operator
$L_{n-1}=h(t,\e)\,\partial^{n-1}+\cdots$ of order $n-1$ with the
leading coefficient $h(\cdot,0)\ne0$, vanishes on the functions
$f_1,\dots,f_{n-1}$, so that $L_{n-1}f_j=0$ for all
$j=1,\dots,n-1$.

By assumption, $L_{n-1}f_n$ is not identically zero, otherwise
$f_n$ would be a linear combination of $f_1,\dots,f_{n-1}$ for all
$\e\ne0$. Therefore, one can represent $L_{n-1}f_n=\e^\nu g(t,\e)$
with some finite $\nu$ and $g(\cdot,0)\ne0$.

Let $L_1$ be the first order operator vanishing on $g$,
constructed as in \eqref{L1}. The composition $L_n=L_1\circ
L_{n-1}$ vanishes on all functions $f_1,\dots,f_{n-1}$ and also on
$f_n$. The coefficients of this operator are analytic functions of
$t,\e$ and the leading term is equal to $g(t,\e)h(t,\e)$ which
does not vanish identically for $\e=0$.
\end{proof}

For the differential operator $L_\e$ constructed in
Lemma~\ref{lode-e}, the ratios
$\|a_j(\cdot,\e)\|/\|a_0(\cdot,\e)\|$ in the sense of \emph{any}
norm on the space of analytic functions on the disk $\mathbb D$,
are automatically bounded as $\e\to0$ since the denominator does
not vanish at $\e=0$. But since $L_\e$ is \emph{unique} modulo a
common factor, this boundedness holds true for \emph{any} analytic
family of operators annulled by the functions $f_1,\dots,f_n$.

More precisely, assume that another differential operator
$L^*_\e=p_0(t,\e)\partial^n+\cdots+p_{n-1}(t,\e)\partial+p_n(t,\e)$
of the same order $n$ has \emph{polynomial} coefficients
$p_j(\cdot,\e)\in\C[t]$ and defines the same set of solutions
$f_1,\dots,f_n$ for all $\e\ne0$. Denote
\begin{equation}\label{bjs}
    b_j^*(\e)=\|p_j(\cdot,\e)\|_{\C[t]}\ge 0.
\end{equation}

\begin{Cor}\label{cor-m}
All the ratios $b_j^*(\e)/b_0^*(\e)$, $j=1,\dots,n$, have finite
limits as $\e\to0$.
\end{Cor}

\begin{proof}
Let $\nu_j\in\mathbb Z$, $j=0,\dots,n$, be the vanishing orders of
the coefficients: this means by definition that
$p_j(t,\e)=\e^{\nu_j}q_j(t,\e)$, where $q_j(\cdot,\e)\in\C[t]$ are
still analytic in $\e$, polynomial in $t$ and $q_j(\cdot,0)\ne0$.
Since any two linear differential operators with the same null
space are proportional, we have
\begin{equation}\label{3rat}
    \frac{p_j}{p_0}=\e^{\nu_j-\nu_0}\frac{q_j}{q_0}
    =\frac{a_j}{a_0},\qquad
    j=1,\dots,n.
\end{equation}
Since the limit values $a_0(\cdot,0)$, $q_j(\cdot,0)$ and
$q_0(\cdot,0)$ are all nonzero by construction, we conclude that
$\nu_j\ge\nu_0$.

By construction,
$b_j^*/b_0^*=\e^{\nu_j-\nu_0}\frac{\|q_j\|_{\C[t]}}{\|q_0\|_{\C[t]}}$,
the second ratio has a bounded limit since the denominator is
nonzero, and the factor $\e^{\nu_j-\nu_0}$ is bounded.
\end{proof}

\subsection{Proof of Lemma~\ref{grig}}
Consider the polynomial family \eqref{meq} satisfied by $n$ entire
functions analytically depending on the parameters, and assume
that some of the ratios $b_j(\l)/b_0(\l)$ is not locally bounded.
Then the semialgebraic set $\{(\l,z)\in\R^{2n^2m}\times\R_+,\
b_0(\l)< zb_j(\l)\}$, a subgraph of the reciprocal function
$b_0(\l)/b_j(\l)$, contains a point $(\l_*,0)$ in its closure. By
the curve selection lemma \cite{milnor:sing-points}, one can find
a real analytic arc $\gamma\:(\R^1,0)\to(\R^{2n^2m},\l_*)$,
$\e\mapsto\l(\e)$, such that along this arc the ratio
$b_0(\l(\e))/b_j(\l(\e))$ tends to zero as $\e\to0$.

But this would contradict the Corollary~\ref{cor-m}, since the
corresponding one-parameter family of linear ordinary differential
equations $L^*_\e=L_{\l(\e)}$ with polynomial coefficients,
obtained by restriction of \eqref{meq} on the arc $\gamma$, has an
analytic family of entire solutions. \qed

\section{Discussion and concluding remarks}

\subsection{Algebra vs.~analysis}
Derivation of the equation \eqref{meq} in \secref{sec:deriv} as a
linear combination between the vectors $\xi_0,\xi_1,\dots,\xi_n$
over the field $\C(t,\l)$ is a distant analog of derivation of a
linear dependence between vectors
$\zeta,B\zeta,B^2\zeta,\dots,B^n\zeta$, where $\zeta\in\C^n$ is a
vector and $B\in\Mat(n,\C)$ a linear operator over the field $\C$.
However, in the latter case the coefficients of the linear
combination are explicitly bounded. Indeed, if $P(u)=u^n+a_1
u^{n-1}+\cdots+a_n$ is the characteristic polynomial of the
operator $B$, then $P(B)=0$ which translates into a linear
dependence between iterates $B^j\zeta$, $j=0,\dots,n$ for
\emph{any} initial vector $\zeta$. The leading coefficient of this
combination is $1$. On the other hand, the magnitude of the
coefficients of the characteristic polynomial is explicitly
bounded in terms of the spectral radius of $B$, i.e., ultimately
in terms of the norm $|B|$. Thus an analog of
Lemma~\ref{cor:bound} in the ``constant'' context, is a purely
algebraic fact.

In a rather surprising way, the bound for the norms of
coefficients in the non-constant case, is obtained using
transcendental arguments, involving solutions of the differential
equations. \emph{A challenging problem is to obtain a direct
algebraic proof of this fact}.

\subsection{Generalizations for Fuchsian systems}
The condition that the system \eqref{ls-p} is polynomial, can be
to some extent relaxed: most of the arguments would work equally
well for \emph{rational} systems exhibiting singularities (regular
or irregular).

Consider, for instance, a Fuchsian system
\begin{equation}\label{fs}
    \dot x=A(t)x,\qquad A(t)=\sum_{j=1}^m\frac{A_j}{t-t_j},
\end{equation}
with the singular points $t_1,\dots,t_m$ and the respective
residue matrices $A_j\in\Mat(n,\C)$ (the point $t=\infty$ is
singular if $\sum_1^m A_j\ne0$). The parameters describing this
system, are $n^m$ entries of the residue matrices, and $m$ complex
numbers $t_1,\dots,t_m$ specifying the location of singular
points. As before, denote by $\l\in\C^{(n^2+1)m}$ the collection
of all these data.

\begin{Thm}\label{fs-thm}
Any linear combination of solutions of the system \eqref{fs}
satisfies a linear $n$th order differential equation of the form
\eqref{lode} with polynomial coefficients $a_j\in\C[t]$,
polynomially depending on the parameters $\l$.

If $|A_j|\le M$, and $|t_j|<M$ for all $j=1,\dots,m$, then
\begin{equation}\label{fs-bd}
    \|a_j\|_{\C[t]}\le\|a_0\|_{\C[t]}\cdot
    [\max(M,2)]^{2^{O(n^2m)}},\qquad \forall j=1,\dots,n.
\end{equation}
\end{Thm}

\begin{proof}[Sketch of the proof]
The construction is completely parallel to the proof of
Corollary~\ref{eq-der} and Lemma~\ref{cor:bound}. One has to
consider the rational vector functions \eqref{iter} in which
$A(t,\l)$ is not polynomial anymore, but rather rational matrix
function whose entries are ratios of lattice polynomials from
$\mathbb Z[t,\l]$. In the same way as before the degrees and the
norms of the numerators and the denominators can be explicitly
controlled. The estimates will be different, but all difference
will be finally absorbed in the $O(n^2m)$-term. To prove the
Principal Lemma~\ref{grig} in this case, one should choose instead
of the unit disk, any other disk free of singularities of the
system \eqref{fs}: for each particular value of the parameters it
is always possible.
\end{proof}

Of course, Lemma~\ref{lem:ode} does not apply to the case when the
\emph{solutions} of the equation have singularities in the unit
disk $\mathbb D$. However, a suitable technical modification of
this Lemma (cf.~with Remark~\ref{rem:larger-disk}) would allow to
place an upper bound for the number of isolated zeros of solutions
of \eqref{lode} in any disk free of singularities of the
solutions, in terms of the size of this disk and the distance from
it to the singular set.

For instance, the following Corollary immediately follows from
Theorem~\ref{fs-thm} and Remark~\ref{rem:larger-disk}.

\begin{Cor}\label{cor:fs}
The number of zeros of any linear combination of solutions of the
system \eqref{fs} in the disk of radius $r<1$ does not exceed the
expression \eqref{main-bd} provided that the concentric disk of
radius $6r$ still contains no singular points of this system. \qed
\end{Cor}

Moreover, the technique developed in \cite{mrl-96} and later in
\cite{quasialg} allows to extend the upper bounds also for
neighborhoods of \emph{regular} singularities with the real
spectrum and sufficiently distant from all \emph{other}
singularities of the system \eqref{fs}. Finally, one could expect
on this way a significant improvement of the bounds given by the
algorithm from \cite{quasialg}.

\subsection{Further applications}
The characteristics that essentially determines the bounds for the
number of zeros in Theorems~\ref{main} and~\ref{fs-thm}, is the
number $n^2m$ (resp., $(n^2+1)m$) of the ``free parameters'' on
which depends the corresponding ``universal family'' of polynomial
(resp., Fuchsian) systems.

Sometimes the ``universal family'' appears naturally in the form
depending on fewer number of parameters. Thus, for instance, the
Picard--Fuchs system of differential equations for Abelian
integrals has a \emph{hypergeometric} form,
\begin{equation}\label{hg}
    (tE-B)\dot x=Cx,\qquad B,C,E\in\Mat(n,\C),\quad E=
    \text{identity matrix},
\end{equation}
completely determined by $2$ constant matrices $B,C$ involving
altogether only $2n^2$ free parameters. The singular points of
this system occur at the eigenvalues of the matrix $B$
\cite{redundant}.

\begin{Thm}\label{hypergeom}
Any linear combination of solutions of the system \eqref{hg}
satisfies a linear $n$th order differential equation of the form
\eqref{lode} with polynomial coefficients $a_j\in\C[t]$,
polynomially depending on the parameters $\l$.

If $|B|,|C|\le M$, then
\begin{equation}\label{fs-hg}
    \|a_j\|_{\C[t]}\le\|a_0\|_{\C[t]}\cdot
    [\max(M,2)]^{2^{O(n^2)}},\qquad \forall j=1,\dots,n.
\end{equation}
\end{Thm}

\begin{proof}
The proof is completely similar to the proof of
Theorems~\ref{main} and~\ref{fs-thm}.
\end{proof}

Note that this is a better bound than the one obtained by the
straightforward reduction of the system \eqref{hg} to the rational
system
\begin{equation*}
    \dot x=\frac{A(t)}{\det(tE-B)}\cdot x,\qquad A(t)\in\Mat(n,\C[t]),
    \ \deg A\le n,
\end{equation*}
involving $n^3+n$ free parameters. Of course, an analog of
Corollary~\ref{cor:fs} is valid in the hypergeometric case as
well.

\subsection{Acknowledgements}
I am greatly indebted to Nicolai Vorobjov for the detailed
explanation concerning effective algorithms for solving polynomial
lattice equations/inequalities. I am also grateful to Dmitry
Novikov with whom we discussed this subject repeatedly, and to
Misha Sodin for valuable comments, especially concerning
Lemma~\ref{lem:ode}.

\bibliographystyle{amsalpha}

\begin{thebibliography}{dlVP29}

\bibitem[BPR03]{roy.e.a.}
S.~Basu, R.~Pollack, and M.-F. Roy, \emph{Algorithms in real
algebraic
  geometry}, Algorithms and Computation in Mathematics, vol.~10,
  Springer-Verlag, Berlin, 2003. \MR{2004g:14064}

\bibitem[dlVP29]{poussin}
C.~de~la Valle{\'e}~Poussin, \emph{Sur l'\'equation
diff\'erentielle lin\'eaire
  du second ordre. d\'etermination d'une int\'egrale par deux valeurs
  assign\'ees. extension aux \'equations d'ordre $n$}, J. Math. Pures Appl.
  \textbf{8} (1929), 125--144.

\bibitem[Gri01]{alexg:thesis}
A.~Grigoriev, \emph{Singular perturbations and zeros of {A}belian
integrals},
  Ph. {D}. thesis, Weizmann Institute of Science (Rehovot), December 2001.

\bibitem[Gri03]{alexg:arxiv}
\bysame, \emph{Uniform asymptotic bound on the number of zeros of
{A}belian
  integrals}, ArXiv preprint \textbf{\arxivno{math.DS/0305248}} (2003).

\bibitem[GV88]{grigoriev-vorobjov}
D.~Yu. Grigor{\cprime}ev and N.~N. Vorobjov, Jr., \emph{Solving
systems of
  polynomial inequalities in subexponential time}, J. Symbolic Comput.
  \textbf{5} (1988), no.~1-2, 37--64. \MR{89h:13001}

\bibitem[HRS90]{heintz-roy-solerno}
J.~Heintz, M.-F. Roy, and P.~Solern{\'o}, \emph{Sur la
complexit\'e du principe
  de {T}arski-{S}eidenberg}, Bull. Soc. Math. France \textbf{118} (1990),
  no.~1, 101--126. \MR{92g:03047}

\bibitem[IY96]{jde-96}
Yu. Ilyashenko and S.~Yakovenko, \emph{Counting real zeros of
analytic
  functions satisfying linear ordinary differential equations}, J. Differential
  Equations \textbf{126} (1996), no.~1, 87--105. \MR{97a:34010}

\bibitem[Lev80]{levin}
B.~Ja. Levin, \emph{Distribution of zeros of entire functions},
revised ed.,
  Translations of Mathematical Monographs, vol.~5, American Mathematical
  Society, Providence, R.I., 1980, Translated from the Russian by R. P. Boas,
  J. M. Danskin, F. M. Goodspeed, J. Korevaar, A. L. Shields and H. P.
  Thielman. \MR{81k:30011}

\bibitem[Mil68]{milnor:sing-points}
J.~Milnor, \emph{Singular points of complex hypersurfaces}, Annals
of
  Mathematics Studies, No. 61, Princeton University Press, Princeton, N.J.,
  1968. \MR{39 \#969}

\bibitem[Nov01]{mit:counterexample}
D.~Novikov, \emph{Systems of linear ordinary differential
equations with
  bounded coefficients may have very oscillating solutions}, Proc. Amer. Math.
  Soc. \textbf{129} (2001), no.~12, 3753--3755 (electronic). \MR{1 860 513}

\bibitem[NY99]{annalif-99}
D.~Novikov and S.~Yakovenko, \emph{Trajectories of polynomial
vector fields and
  ascending chains of polynomial ideals}, Ann. Inst. Fourier (Grenoble)
  \textbf{49} (1999), no.~2, 563--609. \MR{2001h:32054}

\bibitem[NY01]{redundant}
\bysame, \emph{Redundant {P}icard-{F}uchs system for {A}belian
integrals}, J.
  Differential Equations \textbf{177} (2001), no.~2, 267--306. \MR{1 876 646}

\bibitem[NY03]{quasialg}
\bysame, \emph{Quasialgebraicity of {P}icard--{V}essiot fields},
Mosc. Math. J.
  \textbf{3} (2003), no.~2, 551--591, (available as ArXiv preprint
  \arxivno{math.DS/0203210}).

\bibitem[RY96]{mrl-96}
M.~Roitman and S.~Yakovenko, \emph{On the number of zeros of
analytic functions
  in a neighborhood of a {F}uchsian singular point with real spectrum}, Math.
  Res. Lett. \textbf{3} (1996), no.~3, 359--371. \MR{97d:34004}

\bibitem[Yak99]{fields}
S.~Yakovenko, \emph{On functions and curves defined by ordinary
differential
  equations}, The Arnoldfest (Toronto, ON, 1997), Amer. Math. Soc., Providence,
  RI, 1999, pp.~497--525. \MR{2001k:34065}

\bibitem[Yak01]{montreal}
\bysame, \emph{Quantitative theory of ordinary differential
equations and
  tangential {H}ilbert 16th problem}, ArXiv preprint
  \textbf{\arxivno{math.DS/0104140}} (2001), 1--78, Lecture notes of the course
  delivered on the Workshop "Asymptotic series, differential algebra and
  finiteness theorems" (Montreal, June-July, 2000).

\end{thebibliography}
\def\BbbR{$\mathbf R$}\def\BbbC{$\mathbf
  C$}\providecommand\cprime{$'$}\providecommand\mhy{--}\font\cyr=wncyr8\def\Bb%
bR{$\mathbf R$}\def\BbbC{$\mathbf
  C$}\providecommand\cprime{$'$}\providecommand\mhy{--}\font\cyr=wncyr8
\providecommand{\bysame}{\leavevmode\hbox
to3em{\hrulefill}\thinspace}
\providecommand{\MR}{\relax\ifhmode\unskip\space\fi MR }
\providecommand{\MRhref}[2]{%
  \href{http://www.ams.org/mathscinet-getitem?mr=#1}{#2}
} \providecommand{\href}[2]{#2}

\end{document}